\documentclass[11pt]{article}
\usepackage{amssymb,amsfonts,amsthm,ifpdf}
\usepackage{color}
\newtheorem{thm}{Theorem}[section]
\newtheorem{cor}[thm]{Corollary}
\newtheorem{lem}[thm]{Lemma}
\newtheorem{prop}[thm]{Proposition}

\newtheorem{ex}[thm]{Example}

\def\eop{\hfill\rule{2.5mm}{2.5mm}}
\def\pf{\par\smallbreak\noindent {\bf Proof.} \ }

\begin{document}

\title{
{\textbf{\Large{Entropy numbers of Reproducing Hilbert Space of zonal positive definite kernels on compact two-point homogeneous spaces}}\\ } \vspace{-4pt}
\author{
\textsc{Karina Gonzalez}\,
\,\,\&\,\,
\textsc{Thaís Jord\~{a}o}
}}
\date{}

\maketitle 

\begin{center}
\parbox{13 cm}{{\small {\bf Abstract:} We present estimates for the covering numbers of the unit ball of Reproducing Kernel Hilbert Spaces (RKHSs) of functions on  $\mathbb{M}^d$ a $d$-dimensional compact two-point homogeneous space. The RKHS is generated by a continuous zonal/isotropic positive definite kernel. We employ the representation in terms of the Schoenberg/Fourier series expansion for continuous isotropic positive definite kernels, given in terms of a family of orthogonal polynomials on  $\mathbb{M}^d$. The bounds we present carry accurate information about the asymptotic constants depending on the dimension of the manifold and the decay or growth rate of the coefficients of the kernel. The results we present extend the estimates previously known for continuous isotropic positive definite kernels on the $d$-dimensional unit sphere. We present the weak asymptotic equivalence for the order of the growth of covering numbers associated to kernels on  $\mathbb{M}^d$ with a convergent geometric sequence of coefficients. We apply our estimates in order to present a bound for the covering numbers of the spherical Gaussian kernel, and to present bounds for formal examples on $\mathbb{M}^d$.}}
\end{center}

\section{Introduction}
\label{intro}

In this paper we present estimates for the covering numbers of the unit ball of Reproducing Kernel Hilbert Spaces (RKHSs) on a compact two-point homogeneous space $\mathbb{M}^d$ of dimension $d \geq 1$. This releases a follow up paper with the extension of the main results in \cite{GONZALEZ2024128121}, obtained on the unit sphere $\mathbb{S}^d$ of $\mathbb{R}^{d+1}$, for a more general class of manifolds that includes the spherical setting. 

The covering number or Kolmogorov $\epsilon$-entropy (\cite{kolmogorov-MR0112032}) on manifolds are popular concepts that play important role in applied areas such as kernel-based learning algorithms and Gaussian process (\cite{Li-MR1733160,Minh2010-MR2677883,WSS-MR1873936}). The bounds for the covering numbers is a useful and powerful tool in order to estimate the probabilistic error of the statistical nature of the observations from which the algorithms designed are learning from (\cite{smale-MR1864085,Zhou-2003-MR1985575}). Estimates with explicit constants and sharpness included in the asymptotic analysis of finding the bounds for the covering numbers are also very important as illustrated in \cite{kuhn22-MR4474650,ingo-MR4174537}. 

Isotropic positive definite kernels on the spherical setting play important roles in many applications, like kernel methods and spatial statistics (\cite{porcuBerg-MR3619442,FENG2022109287,Montserrat-MR3431424}). The series representation in terms of orthogonal polynomials for isotropic positive definite kernels, known as Schoenberg representation (\cite{Gneiting-MR3102554,Schoenberg-MR0005922}). This representation is widely considered on applications and closely related with the Fourier series expansion of the kernels we work with. The Gaussian Kernel is the main and useful example of spherical isotropic positive kernel (\cite{Minh2010-MR2677883}). 

The framework of the paper is a compact two-point homogeneous space $\mathbb{M}^d$ (\cite{BrownDai2005-MR2119285}). This space is a Riemannian manifold, a compact symmetric space of rank 1, and possess an invariant Riemannian normalized metric $d(\cdot, \cdot)$ such that all geodesics on $\mathbb{M}^d$ are closed and have length $2 \pi$ (\cite{cartan1929determination, Helgason-MR145455}). Let $d\sigma_d$ be the volume element on $\mathbb{M}^d$, we write $L^1(\mathbb{M}^d):=L^1(\mathbb{M}^d, \sigma_d)$ for the Banach space of all measurable functions $f$ on $\mathbb{M}^d$ such that 
$$
\|f\|_1:=\frac{1}{\sigma_d(\mathbb{M}^d)}\int_{\mathbb{M}^d} |f(x)| d\sigma_d x<\infty.
$$
The group motions $SO(\mathbb{M}^d)$ of the manifold $\mathbb{M}^d$ is compact and acts transitively on $\mathbb{M}^d$. Let $p$ be a fixed point (a pole) of the manifold $\mathbb{M}^d$ and consider $SO_p$ be the stationary subgroup of $p$ in $SO(\mathbb{M}^d)$. We denote by $d\alpha$ and $d\alpha_p$ the elements of the normalized Haar measures on the groups $SO(\mathbb{M}^d)$ and $SO_p$, respectively.
In this case, $\mathbb{M}^d=SO(\mathbb{M}^d)/SO_p$ is the orbit space of the compact subgroup $SO_p$ of the orthogonal group $SO(\mathbb{M}^d)$. Due the fact the groups $SO(\mathbb{M}^d)$ and $SO_p$ are compact, these Haar measures are bi-invariant, that is, they are invariant with respect to left and right translations. Moreover, on $\mathbb{M}^d$ the measure $\sigma_d$ is induced by a normalized left Haar measure and the following relation holds
$$
\frac{1}{\sigma_d(\mathbb{M}^d)}\int_{\mathbb{M}^d} f(x) d\sigma_d x=\int_{SO(\mathbb{M}^d)} f(\alpha p) d \alpha, \quad f \in L^1(\mathbb{M}^d).
$$
There are four classes of these spaces and one example. A compact two-point homogeneous space of dimension $d \geq 1$ is the 16-dimensional Cayley's elliptic plane $\mathbb{P}^{16}$ or one of the following manifolds: the $d$-dimensional unit sphere $\mathbb{S}^d$, the $(d+1)$-dimensional real projective space $\mathbb{P}^{d+1}(\mathbb{R})$, the $2d$-dimensional complex projective space $\mathbb{P}^{2d}(\mathbb{C})$, or the $4d$-dimensional quaternion projective space $\mathbb{P}^{4d}(\mathbb{H})$, for $d=1,2,3, \ldots$.

Let $L^2(\mathbb{M}^d):=L^2\left(\mathbb{M}^d, \sigma_d\right)$ be the Hilbert space of square integrable functions $f: \mathbb{M}^d \longrightarrow \mathbb{C}$ with the norm $\|\cdot\|_2$ induced by the inner product
$$
\langle f, g\rangle_2=\int_{\mathbb{M}^d} f(x) \overline{g(x)} d \sigma_d x, \quad f, g \in L^2(\mathbb{M}^d).
$$
The Laplace-Beltrami operator $\mathcal{B}_d$ on $\mathbb{M}^d$ has the spectrum discrete formed by real and non-positive numbers. The radial part $\mathcal{B}$ of the Laplace-Beltrami operator can be written as
$$
\mathcal{B}=(1-x)^{-\alpha}(1+x)^{-\beta} \frac{d}{d x}(1-x)^{1+\alpha}(1+x)^{1+\beta} \frac{d}{d x},
$$
for $\alpha=(d-2)/2$, and $\beta$ is given as follows: for $\mathbb{M}^d=P^{16}$, $\beta=3$, for the case $\mathbb{M}^d=\mathbb{S}^d$ or $\mathbb{M}^d = P^d(\mathbb{R})$, then $\beta = \alpha = (d-2)/2$, if $\mathbb{M}^d=P^d(\mathbb{C})$, then $\beta=0$, and for the remaining case $\mathbb{M}^d=P^d(\mathbb{H})$, we have $\beta=1$. The eigenfunctions of $\mathcal{B}$ are the Jacobi polynomials $P_k^{(\alpha, \beta)}$ with the corresponding eigenvalues $-k(k+\alpha+\beta+1)$, $k=0,1,\ldots$. The eigenvalues of $\Delta=-\mathcal{B}$ are arranged in an increasing order, with each eigenvalue associated to the finite dimensional space $\mathcal{H}_k^d$. The dimension $\tau_k^d:=\dim \mathcal{H}_k^d$ is zero if $\mathbb{M}^d=P^d(\mathbb{R})$ and $k$ is odd, otherwise is given by
\begin{equation}
\tau_k^d =\frac{\Gamma(\beta+1)(2 k+\alpha+\beta+1) \Gamma(k+\alpha+1) \Gamma(k+\alpha+\beta+1)}{\Gamma(\alpha+1) \Gamma(\alpha+\beta+2) \Gamma(k+1) \Gamma(k+\beta+1)}, \quad k=0,1,\ldots.
\label{dim}
\end{equation}
The spaces are mutually orthogonal in $L^2(\mathbb{M}^d)$ and it holds $L^2(\mathbb{M}^d)=\bigoplus_{k=0}^{\infty} \mathcal{H}_k^d$. If we set $\beta_0=\{1\}$, and $\beta_k=\left\{S_{k, j} : j=1,2, \ldots, \tau_k^d\right\}$, for $k=1,2,\ldots$, an orthonormal basis of $\mathcal{H}_k^d$ , then $\{\beta_n\}$ is an orthonormal complete system of $L^2(\mathbb{M}^d)$. This allows to consider the Fourier series expansion of $f\in L^2(\mathbb{M}^d)$ as follows
\begin{equation}
f\sim\sum_{k=0}^{\infty} \sum_{j=1}^{\tau_k^d}\left\langle f, S_{k, j}\right\rangle_2 S_{k, j}.
\label{Fourier.series}
\end{equation}

The relation between the Fourier series expansions of isotropic functions and the so-called Schoenberg series representation is given in terms of the Jacobi polynomials. We present just the necessary facts bellow, more details can be found in \cite[Section 2]{peron-MR4454292}.

For $\alpha, \beta>-1$, we write $P_k^{(\alpha, \beta)}:[-1,1] \longrightarrow \mathbb{R}$ for the Jacobi polynomial of degree $k=0,1,\ldots$ associated to the weight $w(x)=(1-x)^\alpha(1+x)^\beta$. Let $L_1^{\alpha, \beta}:=L_1([-1,1], \omega(t)dt)$ be the space of the measurable functions $f:[-1,1] \longrightarrow \mathbb{R}$ such that
$$
\int_{-1}^1|f(t)|\omega(t) dt<\infty.
$$
The family of orthogonal polynomials $\{P_n^{(\alpha, \beta)}\}$ is a complete orthogonal system on the interval $[-1,1]$ with the weighted measure as bellow
$$
\int_{-1}^1 P_k^{(\alpha, \beta)}(t) P_n^{(\alpha, \beta)}(t)\omega(t) d t=\delta_{k n}h_k^{\alpha,\beta}, \qquad h_k^{\alpha,\beta}=\frac{2^{\alpha+\beta+1}}{(2k+\alpha+\beta+1)} \frac{\Gamma(k+\alpha+1) \Gamma(k+\beta+1)}{\Gamma(k+1) \Gamma(k+\alpha+\beta+1)}, 
$$
$k=1,2,\ldots$ with $\delta_{k n}$ representing the Kronecker delta. If $f \in L_1^{\alpha, \beta}$, then it has the series representation given by
\begin{equation}
f(t)\sim\sum_{k=0}^{\infty} a_k^{\alpha, \beta} \frac{P_k^{(\alpha, \beta)}(t)}{P_k^{(\alpha, \beta)}(1)}, \qquad a_k^{\alpha, \beta}=\frac{P_k^{(\alpha, \beta)}(1)}{h_k^{\alpha, \beta}} \int_{-1}^1 f(t) P_k^{(\alpha, \beta)}(t)(1-t)^\alpha(1+t)^\beta d x,
\label{Fourier-Jacobi.series}
\end{equation}
with
\begin{equation}\label{P(1)}
P_k^{(\alpha, \beta)}(1)=\max _{t \in[-1,1]} P_k^{(\alpha, \beta)}(t)=\frac{\Gamma(k+\alpha+1)}{k! \Gamma(\alpha+1)}, \quad k=0,1,\ldots.
\end{equation}

We will deal with a kernel $K$ on $\mathbb{M}^d$, it means a function $K: \mathbb{M}^d \times \mathbb{M}^d \longrightarrow \mathbb{R}$, that is continuous, isotropic, and positive definite. The kernel $K$ is \textit{zonal} (or \textit{isotropic}), if 
$$
K(x, y)=C(\alpha x, \alpha y), \quad x, y \in \mathbb{M}^d,
$$
all $\alpha \in SO(\mathbb{M}^d)$. It holds that $K$ is an isotropic continuous kernel if it can be written as $K(x, y)=f(\cos (d(x, y)))$, $x, y \in \mathbb{M}^d$, for a continuous function $f:[-1,1] \longrightarrow \mathbb{R}$. This equivalence is carefully proved in \cite{Dunkl-MR2615094} for function on the spherical setting, it is standard the adaptation of the proofs for kernels on a compact two-point homogeneous spaces. The kernel $K$ is \emph{positive definite} (\cite[Chapter 3]{Berg-2012harmonic}) if $K$ is symmetric and for any $n\in\mathbb{N}$, $x_1,\ldots,x_n \in \mathbb{M}^d$ and $c_1, \ldots,c_n\in \mathbb{R}$ it holds 
$$
\sum_{i=1}^{n}\sum_{j=1}^{n}c_{i}c_{j}K(x_{i},x_{j})\geq 0.
$$

The characterization of the Schoenberg type (see \cite[Theorem 1.1]{peron-MR4454292} and \cite{porcuBerg-MR3619442})) asserts that an continuous isotropic positive definite kernel $K$ can be written as a series expansion in terms of the Jacobi polynomial, with summable sequence of non-negative real numbers $\{a_n^{\alpha,\beta}\}$ called the sequence of Schoenberg coefficients. Precisely, let $K: \mathbb{M}^d \times \mathbb{M}^d \longrightarrow \mathbb{R}$ be a kernel such that $K(x, y)=f(\cos (d(x, y)))$, $x, y \in \mathbb{M}^d$, for a continuous function $f:[-1,1] \longrightarrow \mathbb{R}$. The kernel $K$ is positive definite if and only if 
\begin{equation}
K(x, y)=\sum_{k=0}^{\infty} a_k^{\alpha, \beta} \frac{P_k^{(\alpha, \beta)}(\cos (d(x, y)))}{P_k^{(\alpha, \beta)}(1)}, \quad x,y \in \mathbb{M}^ d,    
\label{Schoenberg.series.equation}
\end{equation}
with $\{a_n^{\alpha, \beta}\}$ given in (\ref{Fourier-Jacobi.series}) a summable sequence of non-negative real numbers. 

In this paper $K$ is a continuous isotropic positive definite kernel represented by the series expansion (\ref{Schoenberg.series.equation}).  We write $J_k^{\alpha,\beta}(s)= P_k^{\alpha,\beta}(s)/P_k^{\alpha,\beta}(1)$ for the normalized Jacobi polynomial of degree $k$ associated to the pair $(\alpha,\beta)$, for $k=0,1\ldots$. The Addition Formula, given by the following identity (see \cite[p.405]{BrownDai2005-MR2119285} or \cite[Section 2.3]{Morimoto-MR1641900}, for the spherical setting) 
\begin{equation}
\sum_{j=1}^{\tau_k^d}S_{k,j}(x)S_{k,j}(y)=\tau_k^d J_k^{\alpha, \beta}(\cos(d(x,y))), \quad x, y \in \mathbb{M}^d, \quad k=0,1,\ldots,
\label{addition}
\end{equation}
implies the Fourier series representation of $K$ as follows
\begin{equation}\label{WKlegendre}
K(x,y)=\sum_{k=0}^\infty a_k^{\alpha,\beta} J_k^{\alpha, \beta}(\cos(d(x,y))) = \sum_{k=0}^{\infty} \frac{a_k^{\alpha,\beta}}{\tau_k^d} \sum_{j=1}^{\tau_k^d} S_{k,j} (x) S_{k,j} (y), \quad x,y\in \mathbb{M}^d.
\end{equation}
We assume that $a_k^{\alpha,\beta}\downarrow 0$, as $k\to\infty$. 

For a positive definite kernel $K$ the classical Aronszajn theory of RKHSs (\cite{aron-MR0051437}) asserts that there exists a unique reproducing Hilbert space $(\mathcal{H}_{K},\langle \,\cdot\, , \cdot \, \rangle_K)$ of functions on $\mathbb{M}^d$, satisfying:\\
i. $K(x,\, \cdot\, )\in \mathcal{H}_{K}$ for all $x\in \mathbb{M}^d$;\\
ii. span$\{K(\, \cdot\, ,x) : x\in X\}$  is dense in $\mathcal{H}_{K}$;\\
iii. (Reproducing property) $f(x)=\langle{f},{K(\cdot,x)}\rangle_{K}$, for all $x\in\mathbb{M}^d$ and $f\in \mathcal{H}_{K}$.\\
From the reproducing property and the Cauchy-Schwarz inequality, it is not hard to see that it holds
$$
|f(x)-f(y)|\leq \|f\|_K\|K(\cdot,x)-K(\cdot,y)\|_K=\|f\|_K\left(K(x,x)+K(y,y)-2K(x,y)\right)^{1/2}
$$ 
for any $f\in \mathcal{H}_K$, and for all $x,y\in\mathbb{M}^d$. If $K$ is a continuous kernel, then $\mathcal{H}_{K}\hookrightarrow C(\mathbb{M}^d)$, i.e., we have a continuous embedding of the RKHS $\mathcal{H}_K$ into $C(\mathbb{M}^d)$, the space of continuous functions endowed with the supremum norm $\|\,\cdot \,\|_{\infty}$. In this case, if $K$ is represented as in formula (\ref{WKlegendre}), then the embedding $I_{K}:\mathcal{H}_K \longrightarrow C(\mathbb{M}^d)$ satisfies $\|I_K\|\leq \sum_{k=0}^{\infty}a_k^{\alpha,\beta}$, also due to the reproducing property. Precisely, the interest is on this embedding and the estimates for the related covering numbers. The covering numbers of the operator $I_K$ are defined in terms of the concept of the covering numbers for metric spaces. If $A$ is a subset of a metric space $X$ and $\epsilon>0$, then the \emph{covering numbers} $\mathcal{C} (\epsilon, A)$ are the minimal number of balls of $X$ with radius $\epsilon$ which covers $A$. We denote by $B_K$ the unit ball in $\mathcal{H}_{K}$ and the \emph{covering numbers of the kernel $K$} (or the covering numbers of the RKHS of the kernel $K$) are defined in terms of the covering numbers of the embedding $I_K$, given as follows
$$
\mathcal{C} (\epsilon, I_K):= \mathcal{C}(\epsilon,I_K(B_K)), \quad \epsilon>0.
$$

The content of the paper is organized as follows. In Section \ref{basic}, we present some general facts about the RKHS of isotropic positive definite kernels $K$ represented as in (\ref{WKlegendre}), and general properties for covering numbers. In Section \ref{sec-estimates-rapid-decay}, the bounds for the covering numbers of kernels with the sequence of coefficients possessing rapid decay or decay slower than a geometric progression are presented. As a consequence of the main theorems we present the weak equivalence for order of the growth of  the covering numbers for kernels with a geometric progression as sequence of coefficients, and we obtain bounds for the covering numbers of the Gaussian kernel on the $d$-dimensional unit sphere. In Subsection \ref{subsec-spherical-gaussian}, we present the $d$-dimensional spherical kernel of the Gaussian type, and we apply the estimates obtained previously in order to obtain accurate upper and the lower bounds for the covering numbers of those kernels. In Section \ref{sec-general}, we present the bounds for the covering numbers of kernels with the sequence of coefficients possessing decay or growth comparable to an harmonic sequence. We finish the section with a formal example of kernel on $\mathbb{M}^d$ fulfilling the assumptions of the main theorems, and we present the final comments and the acknowledgments.

\section{Basic theoretical background}\label{basic}

The proofs of the theorems apply general estimates for the covering numbers given em terms of operator norm and for operators with finite dimensional range.

The following proposition presents the characterization of the RKHS in terms of Fourier series expansions, the proof easily follows from the unicity of the RKHS and it is standard in theory of RKHSs. If $K$ is an isotropic positive definite kernel represented by the series expansion (\ref{WKlegendre}), then we will employ the following notation 
$$
\eta_k^d:=\sqrt{\frac{a_k^{\alpha,\beta}}{\tau_k^d}}, \quad k=0,1,\ldots,
$$

\begin{lem}\label{RKHS} If $K$ is an isotropic positive definite kernel represented by the series expansion (\ref{WKlegendre}), then the RKHS $(\mathcal{H}_{K},\langle \,\cdot\, , \cdot \, \rangle_K)$ associated with $K$ is 
$$
\mathcal{H}_K = \left\{g:\mathbb{M}^d \longrightarrow \mathbb{R} :
 g(x)=\sum_{k=0}^{\infty} \eta_k^d \sum_{j=1}^{\tau_k^d} c_{k}^{j} S_{k,j} (x), \quad x\in\mathbb{M}^d, \, \{c_{k}^{j}\} \in \ell^2 , \, j=1,..,\tau_k^d \right\},
$$
endowed the inner product
$$
\langle  g,h \rangle_{K} =\sum_{k=0}^{\infty} \sum_{j=1}^{\tau_k^d} c_{k}^{j} d_{k}^{j}, \quad g,h \in \mathcal{H}_K,
$$
where
$$
g(x)= \sum_{k=0}^{\infty} \eta_k^d \sum_{j=1}^{\tau_k^d} c_{k}^{j} S_{k,j} (x) \quad \mbox{and}\quad h(x)= \sum_{k=0}^{\infty}\eta_k^d \sum_{j=1}^{\tau_k^d} d_{k}^{j} S_{k,j} (x), \quad x\in \mathbb{M}^d.
$$ 
\end{lem}

We consider $P_m$ the orthogonal projections onto the subspace $\mathcal{V}_m$ of $\mathcal{H}_K$ and $P_m^s$ the orthogonal projection onto the subspace and $\mathcal{V}_m^s$, the orthogonal complement of $\mathcal{V}_m$, with the basis
$$\left\{\eta_k^dS_{k,j}: \quad j=1,..,\tau_k^d, \,\, k=0,1,\ldots, m\right\} \quad\mbox{and} \quad \left\{\eta_k^d S_{k,j}: \quad j=1,..,\tau_k^d, \,\, k=m+1, m+2, \ldots\right\},
$$ 
respectively, for each $m=1,2,\ldots$. Then, we can write
$$
f=P_m (f) + P_m^s(f), \quad f\in\mathcal{H}_{K}= \bigoplus_{k=0}^{\infty}\mathcal{V}_k, \quad m=1,2,\ldots,
$$
and we consider the natural embeddings $\mathcal{V}_m \hookrightarrow C(\mathbb{M}^d)$ and $\mathcal{V}_m^s \hookrightarrow C(\mathbb{M}^d)$, keeping the notation $P_m$ and $P_m^s$ for them, respectively. The operator norm of the embeddings, in terms of the Fourier coefficients of the kernel, is given in the following lemma. The proof is standard in the theory and we omitted it here. It can be find in the previous article \cite[Proposition 2.3]{GONZALEZ2024128121} for the spherical setting and the same proof fits here. 

\begin{lem}\label{thmNorms} Let $K$ an isotropic positive definite kernel represented by the series expansion (\ref{WKlegendre}).\ For $m=1,2,\ldots$, the following embeddings
$$
I_K : \mathcal{H}_K \longrightarrow C(\mathbb{M}^d), \quad P_m: \mathcal{V}_m \longrightarrow C(\mathbb{M}^d), \quad \mbox{and} \quad P_m^s: \mathcal{V}_m^s \longrightarrow C(\mathbb{M}^d),
$$
satisfy respectively,
$$
\| I_K \|^2 = \sum_{k=0}^{\infty}a_k, \quad \|P_m\|^2 = \sum_{k=0}^{m}a_k, \quad \mbox{and} \quad
\|P_m^s\|^2 = \sum_{k=m+1}^{\infty}a_k.
$$
\end{lem}

We will employ the notation: for two given functions $f,g: (0,\infty) \longrightarrow \mathbb{R_+}$, $f(\epsilon) \asymp g(\epsilon)$ means the \textit{weak equivalence} and it stands for $f(\epsilon)=O(g(\epsilon))$ and $g(\epsilon)=O(f(\epsilon))$, as $\epsilon\to 0$, we write $f(\epsilon) \approx g(\epsilon)$ for the \textit{strong equivalence}, that is, $f(\epsilon)/g(\epsilon) \to 1$, as $\epsilon\to 0$. 

The Stirling's approximation formula (\cite[p. 24]{artin1964gamma}) 
$$
\Gamma(x)\approx\sqrt{2 \pi} x^{x-1 / 2} e^{-x}, \quad \mbox{as $x\to\infty$},
$$
raises from the identity 
$$
\Gamma(x)=\sqrt{2 \pi} x^{x-1 / 2} e^{-x+\mu(x)}, \quad x>0.
$$
where 
$$
\mu(x)=\sum_{n=0}^{\infty}(x+n+1/2) \log \left(1+\frac{1}{x+n}\right)-1=\frac{\theta}{12 x}, \quad 0<\theta<1.
$$
The formula deliveres the following strong equivalence for the dimension defined in (\ref{dim})
\begin{equation}\label{dimapproxm}
\tau_n^d \approx \frac{2\Gamma(\beta+1)}{\Gamma(\alpha+1) \Gamma(\alpha+\beta+2)} n^{2\alpha+1}=\frac{2\Gamma(\beta+1)}{\Gamma(\alpha+1) \Gamma(\alpha+\beta+2)} n^{d-1}, \quad \mbox{as $n\to\infty$}.
\end{equation}

Immediate consequence of the previous considerations is the following result. 

\begin{prop}\label{dimVm} Let $K$ be an isotropic positive definite kernel represented by the series expansion (\ref{WKlegendre}). If $a_k>0$, for $k=0,1,2,\ldots$, then
$$
\dim \mathcal{V}_m \approx\frac{\Gamma(\beta+1)}{\Gamma(\alpha+2) \Gamma(\alpha+\beta+2)}m^{d}, \quad \mbox{as $m\to\infty$}.
$$
\end{prop}

\pf From identity (\ref{dim}), for $m=1,2,\ldots$, it is clear that
$$
\dim \mathcal{V}_m = \sum_{k=0}^m \tau_k^d = \frac{\Gamma(\beta+1)}{\Gamma(\alpha+1) \Gamma(\alpha+\beta+2) }\sum_{k=0}^m \frac{(2k+\alpha+\beta+1) \Gamma(k+\alpha+1) \Gamma(k+\alpha+\beta+1)}{\Gamma(k+1) \Gamma(k+\beta+1)}, 
$$
and we obtain
$$
\dim \mathcal{V}_m = \frac{\Gamma(\beta+1)}{\Gamma(\alpha+1) \Gamma(\alpha+\beta+2)}\frac{(m+1)(\beta+m+1) \Gamma(m+\alpha+2) \Gamma(m+\alpha+\beta+2)}{(\alpha+1) \Gamma(m+2) \Gamma(m+\beta+2)}.
$$
The Stirling's approximation formula now implies
\begin{eqnarray*}
\dim \mathcal{V}_m \approx \frac{\Gamma(\beta+1)}{e^{2\alpha}\Gamma(\alpha+1) \Gamma(\alpha+\beta+2)(\alpha+1)}b_m, \quad \mbox{as $m\to\infty$},
\end{eqnarray*}
with
\begin{eqnarray*}
b_m = \frac{(m+1)(\beta+m+1) (m+\alpha+2)^{m+\alpha+3/2} (m+\alpha+\beta+2)^{(m+\alpha+\beta+3/2)}}{ (m+2)^{m+3/2} (m+\beta+2)^{m+\beta+3/2}},
\end{eqnarray*}
for $m=1,2,\ldots$. Due the fact that the sequence $\{b_m\}_m$ satisfies $b_m\approx e^{2\alpha}m^{2\alpha+2}$, as $m\to\infty$, we gain
\begin{eqnarray*}
\dim \mathcal{V}_m \approx\frac{\Gamma(\beta+1)}{\Gamma(\alpha+1) \Gamma(\alpha+\beta+2)(\alpha+1)}m^{2\alpha+2}, \quad \mbox{as $m\to\infty$},
\end{eqnarray*}
since $\alpha=(d-2)/2$ and $(\alpha+1)\Gamma(\alpha+1)=\Gamma(\alpha+2)$, the strong equivalence stated follows. \eop


We present the basic and general properties of the covering numbers we are going to use in the following sections. The covering numbers as defined in Section \ref{intro} is the particular case of the following.\ Consider $(X, \|\cdot\|_X), (Y,\|\cdot\|_Y)$ two Banach spaces and $B_X$ the unit ball in $X$, then the \emph{covering numbers} of an operator $T: X \longrightarrow Y$ are given by
$$
\mathcal{C} (\epsilon, T):= \mathcal{C}(\epsilon,T(B_X)), \quad \epsilon>0.
$$
It is not hard to see that if $\| T\| \leq \epsilon $, then $\mathcal{C}(\epsilon, T)=1$, for any $\epsilon>0$, and that $C(\epsilon, T)$ is non-increasing in $\epsilon$.\ The following properties can be found in \cite{Carl-MR1098497,pietsch-MR0917067}. If $S,\, T: X\longrightarrow Y$ and $R: Z \longrightarrow X$ are operators on real Banach spaces, then
\begin{enumerate}
\item[c1.] $\mathcal {C} (\epsilon + \delta, T+S) \leq \mathcal{C} (\epsilon, T) \hspace{0.1cm} \mathcal{C}(\delta, S)$, and
 \item[c2.] $\mathcal{C}(\epsilon \delta,TR)\leq \mathcal{C}(\epsilon, T)\,\mathcal{C}(\delta, R)$, for any $\epsilon, \delta>0$;
\item[c3.] if $n:=\mbox{rank}(T) <\infty$, then $\mathcal{C}(\epsilon, T)\leq \left( 1+2\| T\|/\epsilon \right)^ n$, for any $\epsilon>0$.
\end{enumerate}
For the lower estimates we will need the following inequality for the covering numbers for two $n$-dimensional Hilbert spaces $X$ and $Y$, 
\begin{equation}\label{CN1}
 \sqrt{\det{T^{*}T}} \left(\frac{1}{\epsilon}\right) ^n \leq \mathcal{C}(\epsilon,T), \quad \epsilon>0.
\end{equation}

\section{Estimates of kernels with rapid decay of the coefficients}
\label{sec-estimates-rapid-decay}

In this section present the bounds for the covering numbers of the RKHS of kernels on $\mathbb{M^d}$ represented by the series expansion (\ref{WKlegendre}) with sequence of coefficients with decay comparable to a geometric progression. We apply them in order to estimate the covering numbers of the Gaussian kernel on the unit sphere and the covering numbers of spherical kernels the Gaussian type.

In order to present the results we fix the following short notation for the quantities given in Lemma \ref{thmNorms},
\begin{equation}\label{kappas}
\kappa= \|I_K\|, \quad  \kappa_m= \|P_m\|, \quad \mbox{and} \quad \kappa_m^s= \|P_m^s \|, \quad \mbox{for $m=1,2,\ldots$}.
\end{equation}

Let $K$ be a  kernel represented by the series expansion (\ref{WKlegendre}) on $\mathbb{M}^d$, we will write $a_k:=a_k^{\alpha,\beta}$, $k=0,1,\ldots$, for the sequence of coefficients of the kernel. In this case, we observe that if $a_0\neq 0$ and there exists $0<\theta<1$ such that  $a_k \leq \theta\, a_{k-1}$, for $k=1,2\ldots$, then  
\begin{equation}\label{CS}
0<(k_m^s)^2 \leq \frac{\theta^{m+1}}{1-\theta} \, a_0, \quad m=1,2, \ldots.
\end{equation} 

The next theorem is the extension of the Proposition 3.1 of the reference \cite{GONZALEZ2024128121}. 

\begin{thm}\label{thmuppergeo}  Let $K$ a continuous isotropic positive definite kernel given by the series expansion (\ref{WKlegendre}). If there exists $0<\theta<1$ such that  $a_k \leq \theta\, a_{k-1}$, for $k=1,2\ldots$, then 
$$
\limsup_{\epsilon\to 0^+}\frac{\ln (\mathcal{C}(\epsilon, I_K ))}{[\ln(1/\epsilon)]^{d+1}}\leq\frac{2^{d+1}\Gamma(\beta+1)}{\Gamma(\alpha+2) \Gamma(\alpha+\beta+2)}\frac{1}{ [\ln({1}/{\theta})]^d} . 
$$
\end{thm}

\pf 
With no loss of generality we can consider $a_0\neq 0$.\ Otherwise, in the proof bellow, we just replace $a_0$ by the first positive term of the sequence $\{a_n\}$.

For $\epsilon >0$, we consider $m:=m(\epsilon)\in \mathbb{Z}_+$ as follows
\begin{equation}\label{alphaMEpsilonexp}
     \left(\frac{\theta^{m+1}}{1-\theta} \,a_0\right)^{1/2}    < \frac{\epsilon}{2} < \left(\frac{\theta^{m}}{1-\theta}\, a_0\right)^{1/2}.
\end{equation} 
Applying the logarithm above, we have that 
\begin{equation}\label{m1}
m\approx \frac{\ln (\epsilon^2 (1-\theta) / 4\, a_0 )}{\ln\theta}= \frac{2\ln(2a_0^{1/2}/\epsilon(1-\theta)^{1/2})}{\ln(1/\theta)},  \quad \mbox{as $\epsilon\to 0^+$}.  
\end{equation}
Considerations above, inequality (\ref{CS}), and inequality (\ref{alphaMEpsilonexp}) imply 
$$
\mathcal{C}\left(\epsilon , P_m^s \right)=1, \quad m=1,2,\ldots. 
$$
Thus, properties c1. and c3. for covering numbers, respectively, lead us to the estimate
$$
\mathcal{C} (\epsilon, I_K = P_m+P_m^s) \leq \left(1+ 4\| P_m \|/\epsilon \right)^ {\mbox{rank}(P_m)} \leq \left(1+ 4\kappa_m/\epsilon \right)^ {\dim{\mathcal{V}_m}},
$$
that implies
\begin{eqnarray*}
    \ln (\mathcal{C}(\epsilon, I_K ))\leq \dim{\mathcal{V}_m} \ln \left(8\kappa/\epsilon\right), \quad m=1,2,\ldots,
\end{eqnarray*}    
since $\kappa_m\leq\kappa$, for any $m$. From Proposition \ref{dimVm}, we have 
$$
\mathcal{C} (\epsilon, I_K) \leq \dim{\mathcal{V}_m}\ln \left(8\kappa/\epsilon\right)< 2\frac{\Gamma(\beta+1)}{\Gamma(\alpha+1) \Gamma(\alpha+\beta+2)(\alpha+1)}m^{d} \ln \left(8\kappa/\epsilon\right),
$$
for $m=m_0, m_0+1, \ldots$, and some $m_0\in\mathbb{N}$. With the application of the strong equivalences (\ref{m1}) and $\ln t\approx\ln 8\kappa t$, we obtain 
\begin{eqnarray*}
\limsup_{\epsilon\to 0^+}\frac{\ln (\mathcal{C}(\epsilon, I_K ))}{[\ln(1/\epsilon)]^{d+1}}=\limsup_{\epsilon\to 0^+}\frac{\ln (\mathcal{C}(\epsilon, I_K ))}{[\ln(2a_0^{1/2}/\epsilon(1-\theta)^{1/2})]^{d+1}} & \leq & \frac{2^{d+1}\Gamma(\beta+1)}{\Gamma(\alpha+2) \Gamma(\alpha+\beta+2)}\frac{1}{[\ln(1/\theta)]^d}.
\end{eqnarray*}
\eop

We present an application of Theorem \ref{thmuppergeo} to the spherical Gaussian Kernel. The background material we will employ can be found in the references \cite{Jordao-2019-MR3923505,Minh2010-MR2677883}.

\begin{ex}
 The Gaussian kernel on the unit sphere $\mathbb{S}^d$, for $d\geq 2$, is given by
\begin{equation}\label{gaussian-def}
K_\varrho(x, y)=\exp \left(-2 \varrho^{-2}(1-x \cdot y)\right), \quad x, y \in \mathbb{S}^d,    
\end{equation}
for $\varrho>0$, where $x \cdot y$ is the dot product in $\mathbb{R}^{d+1}$. It is well known (see \cite[Theorem 8]{Minh2010-MR2677883} and \cite[Example 2.5]{Jordao-2019-MR3923505}) that the Gaussian kernel is an isotropic (continous) positive definite kernel represented by the series expansion (\ref{WKlegendre}), as follows
\begin{equation}\label{WKlegendrespherical}
K_{\varrho}(x,y)=\sum_{k=0}^{\infty} \lambda_k^\varrho \tau^d_k G_k^{(d-2)/2}(x\cdot y) = \sum_{k=0}^{\infty}\lambda_k^\varrho \sum_{j=1}^{\tau^d_k} S_{k,j} (x) S_{k,j} (y), \quad x,y\in \mathbb{S}^d,
\end{equation}
where $G_k^{(d-2)/2}(s)= J_k^{\alpha,\alpha}(s)$, for $\alpha=(d-2)/2$, is the normalized Gegenbauer polynomial of degree $k$ associated to the dimension $d$, and
\begin{equation}\label{GaussianCoef}
\lambda_k^\varrho=e^{-2 / \varrho^2} \varrho^{d-1} \Gamma\left(\frac{d+1}{2}\right) I_{k+(d-1) / 2}\left(\frac{2}{\varrho^{2}}\right), \quad k=0,1,2, \ldots,    
\end{equation}
with $I_{\nu}(\cdot)$ denoting the modified Bessel function of first kind associated with the parameter $\nu$. The modified Bessel function of first kind is defined in terms of the Bessel function $J_\nu$ with the parameter $\nu>-1/2$ (\cite[p. 5, formula (12)]{ErdelyiII-MR0698780}) and given by
$$
 I_\nu(z)=e^{-i\nu \pi/2} J_\nu (z e^{i \pi/ 2})=\sum_{j=0}^{\infty}\left(\frac{z}{2}\right)^{\nu+2 j} \frac{1}{j! \Gamma(j+\nu+1)}.
$$   
The definition of the Bessel function above implies the following
\begin{eqnarray*}
I_{k+(d-1)/2}\left(\frac{2}{\varrho^2}\right) 
< \frac{1}{\varrho^2(k+(d-1)/ 2)} I_{k-1+(d-1)/2}\left(\frac{2}{\varrho^2}\right),\quad k=0,1,2, \ldots,
\end{eqnarray*}
this means that 
\begin{equation}\label{aaa}
\varrho^2 (k+(d+1)/2) \lambda_{k+1}^\varrho <\lambda_k^\varrho, \quad k=0,1, \ldots.    
\end{equation}
\end{ex}

\begin{cor}
    Let $K_{\varrho}$ be the Gaussian kernel on $\mathbb{S}^d$ defined in (\ref{gaussian-def}). If $\varrho^2>2$, then 
    $$
\limsup _{\epsilon \rightarrow 0^{+}} \frac{\ln \left(\mathcal{C}\left(\epsilon, I_K\right)\right)}{[\ln (1 / \epsilon)]^{d+1}}\leq \frac{2^{d+2}}{d!\,[{\ln\left(\varrho^2 /2\right)}]^d} = \frac{4}{d!} \frac{1}{[\ln(\varrho/ \sqrt{2} )]^d}.
$$
\end{cor}

\pf In this spherical case,  ($\alpha=\beta=(d-2)/2$ in formula (\ref{dim})) we have that 
$$
\tau_k^d =\frac{(2k+d-1)(k+d-2)!}{k!(d-1)!}, 
$$
then the identity
\begin{equation}\label{tau-k}
\tau_{k+1}^d = \frac{(2k+d+1)(k+d-1)}{(2k+d-1)(k+1)}\tau_k^d, \quad \quad k=0,1,2,\ldots,   
\end{equation}
and inequality (\ref{aaa}), imply
\begin{eqnarray*}
\lambda_k^\varrho \tau_k^d  >  \varrho^2 (k+(d+1)/2) \lambda_{k+1}^\varrho\tau_k^d & = & \frac{\varrho^2}{2} \frac{(2k+d-1)(k+1)}{(k+d-1)} \lambda_{k+1}^\varrho \tau_{k+1}^d,
\end{eqnarray*}
for $k=0,1,2,\ldots$. Therefore, it holds the following
\begin{equation}\label{bbb}
 \lambda_{k+1}^\varrho \tau_{k+1}^d <  \frac{2}{\varrho^2}\frac{(k+d-1)}{(2k+d-1)(k+1)} \lambda_k^\varrho \tau_k^d \quad \mbox{for $k=0,1,2,\ldots$}.
\end{equation}

It is easy to see that $\{a_k=\lambda_k^{\varrho}\tau_k^d\}$ is the sequence of coefficients of the kernel $K_{\varrho}$ in the series representation given in formula (\ref{WKlegendre}). Then, the inequality (\ref{bbb}) promptly implies that the sequence of coefficients of the kernel $K_{\varrho}$ satisfies $a_{k+1}<(2/\rho^2) a_k$, $k=1,2,\ldots$. An application of Theorem \ref{thmuppergeo} with $\alpha=\beta=(d-2)/2$ and for $\theta=2/\rho^2$ leads to the estimate stated. \eop

The lower estimate for the covering numbers of kernels on $\mathbb{M}^d$ with decay not faster than a rapid decay is is presented in the next theorem. This extends Proposition 3.2 of \cite{GONZALEZ2024128121}. 

\begin{thm}\label{thmlowergeo}  Let $K$ a continuous isotropic positive definite kernel given by the series expansion (\ref{WKlegendre}). If there exists $0<\delta<1$ such that  $a_k \geq \delta\, a_{k-1}>0$, for $k=1,2\ldots$, then 
\[
 \frac{\Gamma(\beta+1)}{\Gamma(\alpha+2) \Gamma(\alpha+\beta+2)}\,\frac{ 2^{d-1} d^d }{ [\ln(1/\delta)+ d-1]^d (d+1)^{d+1}} \leq \liminf_{\epsilon\to 0^+}\frac{\ln(\mathcal{C}(\epsilon, I_K ))}{[\ln(1/\epsilon)]^{d+1}} .
\]
\end{thm}

\pf For $m=1,2,\ldots$ we consider the following composition operator
$$
T_m : \mathcal{V}_m \stackrel{ J_m }{\longrightarrow} \mathcal{H}_K  \stackrel{I_K}{\longrightarrow} C(\mathbb{M}^d) \stackrel{{L_m}}{\longrightarrow}  L^2 ( \mathbb{M}^d ) \stackrel{P_m }{\longrightarrow} L_m\circ I_K\circ J_m(\mathcal{V}_m),
$$
where $J_m$ stands for the identity operator given by $\mathcal{V}_m\hookrightarrow\mathcal{H}_K$, $L_m$ and
$P_m $ are the orthogonal projections on $L_m I_K J_m(\mathcal{V}_m)$.\ It is clear that $\|P_m\|=\|J_m\|=\Vert L_m \Vert=1$, and 
\begin{equation}\label{det}
\det(T_m ^* T_m)= \prod_{k=0}^{m} \left(\frac{a_k}{\tau_{k}^{d}}\right)^{\tau^{d}_{k}}.
\end{equation}
For any $\epsilon>0$, properties c1. and c2, for the covering numbers, imply 
\begin{eqnarray*}
\mathcal{C} (\,\epsilon, T_m ) =\mathcal{C} (\,\epsilon, P_m L_m I_K J_m ) & \leq & \mathcal{C} (1, P_m ) \mathcal{C} (1,L_m ) \mathcal{C} (\epsilon, I_K )\mathcal{C} (1, J_ m ) = \mathcal{C} (\epsilon, I_K ),
\end{eqnarray*}
for $m=1,2,\ldots$. The application of the lower bound for the covering number given in inequality (\ref{CN1}), leads us to
$$
\sqrt{ \det (T_m ^* T_m )}\left(\frac{1}{{}\,\epsilon}\right)^{\dim \mathcal{V}_m} \leq \mathcal{C} (\,\epsilon, T_m  ) \leq \mathcal{C} (\epsilon, I_K ), 
$$
and due the formula (\ref{det}), we have
$$
\prod_{k=0}^{m} \left(\frac{a_k}{\tau^d_k}\right)^{\tau^d_k/2}\left(\frac{1}{\,\epsilon}\right)^{\dim \mathcal{V}_m}\leq \mathcal{C} (\epsilon, I_K ), \quad m=1,2,\ldots.
$$

For $\epsilon >0$, consider
\begin{equation}\label{J_m}
J_m:= \frac{1}{2}\sum_{k=0}^{m} \tau_k^d \ln \left(a_k/\tau_k^d \right) - \dim\mathcal{V}_m \ln({\epsilon})\leq \ln(\mathcal{C}(\epsilon, I_K )), \quad m=1,2,\ldots,
\end{equation}
We observe that, for $m=1,2,\ldots$, the estimate 
$$
\delta^{m}a_0/\tau_{m}^{d}\leq a_m/\tau_{m}^{d}\leq a_k/\tau_k^{d} , \quad k=0,1,\ldots m,
$$
implies 
\begin{equation}\label{lower2geo}
\frac{ \sum_{k=0}^{m} \tau_{k}^{d} }{2} \left( \ln (\delta^{m}a_0)-\ln (\tau_m^d )\right) \leq \frac{1}{2}\sum_{k=0}^{m} \tau_d^{k} \ln\left(a_k/\tau_k^{d} \right), \quad m=1,2,\ldots.
\end{equation}
If we write
$$
I_m:=\frac{ \sum_{k=0}^{m} \tau_{k}^{d} }{2} \left( \ln (\delta^{m}a_0)-\ln (\tau_m^d)\right) , \quad m=1,2,\ldots,
$$
then,
\begin{eqnarray*}
J_m =I_m - \dim V_m \ln(\epsilon) \leq \ln(\mathcal{C} (\epsilon, I_K)), \quad m=1, 2,\ldots. 
\end{eqnarray*}
Formula (\ref{dimapproxm}) implies that
$$
\tau^d_k \leq \tau^d_m < \frac{3\Gamma(\beta+1)}{\Gamma(\alpha+1) \Gamma(\alpha+\beta+2)}m^{d-1}, 
$$
for $m=m_0, m_0+1,\ldots$ and some $m_0\in\mathbb{Z}_+$. By an application of the Proposition \ref{dimVm}, we can consider that  
$$
\frac{\Gamma (\beta+1)}{2\Gamma(\alpha +2)\Gamma(\alpha+\beta+2)} m^d < \dim V_m = \sum_{k=0}^{m} \tau_k^d , \quad m=m_0, m_0+1,\ldots.
$$
Therefore, for any $\epsilon>0$ and $m=m_0,m_0+1,\ldots$, it holds
$$
A_{\alpha,\beta}\,m^d \left[m\ln \delta + \ln a_0 - (d-1)m - \ln \left(\frac{3\Gamma(\beta+1)}{\Gamma(\alpha+1) \Gamma(\alpha+\beta+2)}\right) - 2\ln(\epsilon)\right]< J_m,
$$
with
$$
A_{\alpha, \beta}:= \frac{\Gamma(\beta+1)}{4\Gamma(\alpha+2) \Gamma(\alpha+\beta+2)}.
$$


We consider 
$$ 
\varphi_\epsilon (m):= -m^d  \ln{\left(\frac{3\Gamma(\beta+1)\,\epsilon^2 }{a_0 \Gamma(\alpha+1) \Gamma(\alpha+\beta+2)}\right)}-m^{d+1}\left[\ln(1/\delta)+ d-1\right],
$$
with $\epsilon>0$ and $m=m(\epsilon)$ such that $1/(m+1)\leq\epsilon <1/m$.\\
The final steps of the proof is to solve the optimization problem given by $\varphi_\epsilon$. Thereby, we obtain
\begin{equation}\label{lower5geo}
A_{\alpha, \beta} \varphi_\epsilon (m) < \ln(\mathcal{C}(\epsilon, I_K )), \quad m=m_0, m_0+1, \ldots,
\end{equation}
for some $m_0\in\mathbb{Z}_+$.
The critical point of $\varphi_{\epsilon}$, as a function on $\mathbb{R}$, is given by
$$
c=-\frac{d}{(d+1)[\ln(1/\delta)+d-1]} \ln\left(\frac{3\Gamma(\beta+1)\epsilon^2}{a_0\Gamma(\alpha+\beta+2)\Gamma(\alpha+1)}\right).
$$
Since $\lfloor c \rfloor \approx  c \approx \lceil c \rceil $, we can estimate (\ref{lower5geo}) in terms of $\varphi_{\epsilon}(c)$, even if $c$ is not an integer.\ A simple calculation leads us to 
\begin{eqnarray*}
A_{\alpha,\beta} \varphi_\epsilon (c) =A_{\alpha, \beta} \frac{[- \ln(3 \,\Gamma(\beta+1)\,\epsilon^2 /a_0 \Gamma(\alpha+1) \Gamma(\alpha+\beta+2)]^{d+1}}{ [\ln(1/\delta)+ d-1]^d} \frac{d^d}{(d+1)^{d+1}},
 \end{eqnarray*}
we can consider $a_0$ such that $a_0 \Gamma(\alpha+1) \Gamma(\alpha+\beta+2)/3\, \Gamma(\beta+1)\geq 1$, and then inequality (\ref{lower5geo}) implies
$$
\frac{ 2^{d+1} d^d A_{\alpha,\beta}}{[\ln(1/\delta)+ d-1]^d (d+1)^{d+1}}\left[\ln\left(\frac{\sqrt{a_0\Gamma(\alpha+\beta+2)\Gamma(\alpha+1)}}{\sqrt{3\Gamma(\beta+1)}\epsilon}\right)\right]^{d+1} < \ln(\mathcal{C}(\epsilon, I_K )),
$$
for $1/(m+1)\leq\epsilon <1/m$ and $m=m_0, m_0+1, \ldots$. The proof follows by observing that 
$$
\ln\left(\frac{\sqrt{a_0\Gamma(\alpha+\beta+2)\Gamma(\alpha+1)}}{\sqrt{3\Gamma(\beta+1)}\epsilon}\right)\approx \ln(1/\epsilon), \quad \mbox{as $\epsilon\to 0^+$.}
$$
\eop

The next consequence of the Theorem \ref{thmuppergeo} and the Theorem \ref{thmlowergeo} is the weak equivalence for the asymptotic behavior of the covering numbers of kernels with a convergent geometric sequence of coeffcients. This equivalence is well known for the spherical kernels with a convergent geometric sequence of coeffcients (\cite[Theorem 1.1]{GONZALEZ2024128121}) and it deliveres us the exact order of the growth of the covering numbers. 

\begin{cor}\label{corweakequiv}   Let $K$ a continuous isotropic positive definite kernel given by the series expansion (\ref{WKlegendre}), with $a_k=\delta^k$, for $k=1,2\ldots$, for some $0<\delta<1$. Then, 
\[
\ln(\mathcal{C}(\epsilon, I_K )) \asymp [\ln(1/\epsilon)]^{d+1}, \quad \mbox{as $\epsilon\to 0^+$}.
\]
\end{cor}

\subsection{Estimates of spherical kernels of the Gaussian type} \label{subsec-spherical-gaussian}

We close this section with the application of the Theorem \ref{thmuppergeo} and the Theorem \ref{thmlowergeo} to estimate the covering numbers of the RKHS associated to the spherical kernel of the Gaussian type.

A kernel of the Gaussian type on the unit sphere $\mathbb{S}^d$, for $d\geq 1$, is defined by
\begin{equation}\label{WKlegendrespherical-like}
K_{\delta}(x,y)=\sum_{k=0}^{\infty} \lambda_k \tau^d_k G_k^{(d-2)/2}(x\cdot y) = \sum_{k=0}^{\infty}\lambda_k \sum_{j=1}^{\tau^d_k} S_{k,j} (x) S_{k,j} (y), \quad x,y\in \mathbb{S}^d,
\end{equation}
with $\{\lambda_k:=\delta^k\}$ for some $0<\delta<1$. 

Concerning about the estimates for the covering numbers of the RKHS of kernels of the Gaussian type we have the following upper and lower estimates.

\begin{prop}\label{propGaussianespherical}
Let $K_{\delta}$ be a kernel of the Gaussian type represented by the series expansion (\ref{WKlegendrespherical-like}) on $\mathbb{S}^d$. For any $0<\delta<1/(d+3)$, if $d\geq 2$, and for any $0<\delta<1$, if $d=1$, it holds
$$
\frac{1}{2(d+1)!}\left[\frac{d/(d+1)}{\ln(e^{d-1}/\delta)}\right]^d \leq \liminf_{\epsilon\to 0^+}\frac{\ln(\mathcal{C}(\epsilon, I_K ))}{[2\ln(1/\epsilon)]^{d+1}} \leq \limsup _{\epsilon \rightarrow 0^{+}} \frac{\ln \left(\mathcal{C}\left(\epsilon, I_K\right)\right)}{[2\ln (1 / \epsilon)]^{d+1}}\leq \frac{2}{ d!}\frac{1}{[{\ln\left( {1}/{\delta(d+3)}\right)}]^d}.
$$
\end{prop}

\pf For the case $d\geq 2$, the sequence of coefficients $\{a_k=\delta^k\tau_k^d\}$ fits in the assumptions of Theorem \ref{thmuppergeo} and Theorem \ref{thmlowergeo}, for the given $0<\delta<1/(d+3)$ and $\theta=\delta(d+3)$, respectively. In fact, the application of identity (\ref{tau-k}) leads us to
\begin{eqnarray*}
a_{k+1}=\delta^{k+1}\tau^d_{k+1}=\delta\frac{(2k+d+1)(k+d-1)}{(2k+d-1)(k+1)}a_k, \quad k=0,1,\ldots.  
\end{eqnarray*}
The inequality
$$
1\leq \frac{(2k+d+1)(k+d-1)}{(2k+d-1)(k+1)}\leq \frac{2k+d+1}{k+1}\leq d+3, \quad k=0,1,\ldots,
$$
implies that
\begin{eqnarray}\label{gaussian-like-d}
\delta a_k\leq a_{k+1}\leq \delta (d+3)a_k, \quad k=0,1,\ldots.  
\end{eqnarray}
We remember that for $\mathbb{M}^d=\mathbb{S}^d$, $d\geq 1$, we have $\alpha=\beta=(d-2)/2$. Then, we obtain
$$
\frac{ 2^{d} d^{d} }{ [\ln(1/\delta)+ d-1]^d d!(d+1)^{d+1}} \leq \liminf_{\epsilon\to 0^+}\frac{\ln(\mathcal{C}(\epsilon, I_K ))}{[\ln(1/\epsilon)]^{d+1}},
$$
and
$$
\limsup _{\epsilon \rightarrow 0^{+}} \frac{\ln \left(\mathcal{C}\left(\epsilon, I_K\right)\right)}{[\ln (1 / \epsilon)]^{d+1}}\leq \frac{2^{d+2}}{ d!}\frac{1}{[{\ln\left(1/\delta(d+3)\right)}]^d}.
$$
A simple rearrangement of the terms in the inequalities above leads us to the inequality stated in the proposition. 

For $d=1$ the sequence of coefficients $\{a_k=\delta^k\tau_k^1\}$ fits in the assumptions of Theorem \ref{thmuppergeo} and Theorem \ref{thmlowergeo}, for $\theta=\delta$. In fact, we know that $\tau_k^1=2$, for $k=1,2,\ldots$, then 
\begin{eqnarray}\label{gaussian-like-1}
a_{k+1}=\delta^{k+1}\tau^1_{k+1}=2\delta^{k+1}=\delta \delta^k\tau_k^1=\delta a_k, \quad k=0,1,\ldots.
\end{eqnarray}
We remember that $\alpha=\beta=-1/2$, then we obtain
$$
\frac{1}{2}\,\frac{ 1 }{ \ln(1/\theta)} \leq \liminf_{\epsilon\to 0^+}\frac{\ln(\mathcal{C}(\epsilon, I_K ))}{[\ln(1/\epsilon)]^{2}} \leq \limsup _{\epsilon \rightarrow 0^{+}} \frac{\ln \left(\mathcal{C}\left(\epsilon, I_K\right)\right)}{[\ln (1 / \epsilon)]^{2}}\leq \frac{8}{{\ln\left( {1}/{\theta}\right)}}\leq \frac{8}{{\ln\left( {1}/{4\theta}\right)}}.
$$
This is the inequality stated. \eop

\section{Estimates of kernels with general decay of the coefficients}
\label{sec-general}

In this section present the bounds for the covering numbers of kernels on $\mathbb{M}^d$, represented by the series expansion (\ref{WKlegendre}), with sequence of coefficients with decay comparable to an harmonic progression.

We observe that the kernel $K$ is represented by the series (\ref{WKlegendre}) with sequence of coefficients $\{a_k\}$ satisfying
$$\limsup_{n\to\infty}a_n/n^{-\gamma}\leq c,$$ for some $\gamma>1$ and $c>0$, then there exists $m_0\in\mathbb{N}$ such that
\begin{equation}\label{CSpol}
0< (k_m^s)^2 \leq \frac{c}{\gamma-1}{(m+1)}^{-\gamma+1}, \quad m=m_0,m_0+1, \ldots,
\end{equation} 
with $\kappa_m^s$ as introduced in (\ref{kappas}). Supposing only $K$ given as (\ref{WKlegendre}) on $\mathbb{M}^d$ with summable sequence of coefficients $\{a_k\}$, then $a_k/\tau_k^d = O(k^{-d})$ (as on the spherical setting \cite[Section 5]{GONZALEZ2024128121}) and due the formula (\ref{dimapproxm}), it holds that $a_k=O(k^{-1})$. We consider this basic decay for our assumption in the next result. 

The next theorem extends Proposition 3.3 of the reference \cite{GONZALEZ2024128121}. The proofs are similar to the previous case (Section \ref{sec-estimates-rapid-decay}) and we will short the presentation for the similar details.

\begin{thm}\label{thmupperhar} Let $K$ a continuous isotropic positive definite kernel given by the series expansion (\ref{WKlegendre}).\ If there exists $\gamma>0$ and $c_1>0$ such that $a_k/ k^{-d-\gamma} \leq c_1$, for $k=n_0, n_0+1, \ldots$, then 
$$
\limsup_{\epsilon\to 0^+}\frac{\ln (\mathcal{C} (\epsilon, I_K))}{(1/\epsilon)^{2d/(\gamma+d-1)}\ln (1/\epsilon)}\leq \frac{\Gamma(\beta+1)}{\Gamma(\alpha+2) \Gamma(\alpha+\beta+2)} \left(\frac{4 c_1 }{\gamma+d-1}\right)^{d/(\gamma+d-1)}. 
$$
\end{thm}
\pf For $\epsilon>0$, we consider $m:=m(\epsilon)\in \mathbb{N}$ such that
\begin{equation}\label{alphaMEpsilonpol}
     \left[\frac{c_1}{\gamma+d-1}(m+1)^{-\gamma-d+1}\right]^{1/2}    < \frac{\epsilon}{2 }< \left[\frac{c_1}{\gamma+d-1}m^{-\gamma-d+1}\right]^{1/2}.
\end{equation} 
In this case, we have
\begin{equation}\label{mEpsilonpol}
m(\epsilon) \approx \left(\frac{4 c_1 }{\gamma+d-1}\right)^{1/(\gamma+d-1)} \left( \frac{1}{\epsilon}\right)^{2/(\gamma+d-1)}, \quad \mbox{as $\epsilon\to 0^+.$}
\end{equation}
Repeating exactly the same steps in the proof of Theorem \ref{thmuppergeo}, we obtain
$$
\mathcal{C} (\epsilon, I_K) \leq \dim{\mathcal{V}_m}\ln \left(8\kappa/\epsilon\right)\approx \frac{\Gamma(\beta+1)}{\Gamma(\alpha+2) \Gamma(\alpha+\beta+2)}m^{d} \ln \left(8\kappa/\epsilon\right),
$$
as $m\to\infty$. With the application of the strong equivalences (\ref{mEpsilonpol}), we obtain 
$$
\limsup_{\epsilon\to 0^+}\frac{\ln (\mathcal{C}(\epsilon, I_K ))}{(1/\epsilon)^{2d/(\gamma+d-1)}\ln (1/\epsilon)}  \leq \frac{\Gamma(\beta+1)}{\Gamma(\alpha+2) \Gamma(\alpha+\beta+2)} \left(\frac{4\,c_1 }{\gamma+d-1}\right)^{d/(\gamma+d-1)}.
$$ 
\eop

The next result is the extension of the Proposition 3.4 of the reference \cite{GONZALEZ2024128121}.

\begin{thm}\label{thmlowerhar} Let $K$ a continuous isotropic positive definite kernel given by the series expansion (\ref{WKlegendre}).\ If there exist $\rho>1$ and $c_2>0$ such that $a_k/k^{-\rho} \geq c_2$, for $k= 1,2,\ldots$, then
$$ 
\ln{\sqrt{a_0}}+ \frac{e^{-1}(\rho+d-1)\Gamma (\beta+1)}{d \,\Gamma(\alpha+2) \Gamma(\alpha+\beta +2) }\left( \frac{c_2 \Gamma(\alpha+1) \Gamma(\alpha+\beta+2)}{3\Gamma(\beta+1)} \right)^{d/(\rho+d-1)}\leq \liminf_{\epsilon\to 0}\frac{\ln(\mathcal{C}(\epsilon, I_K ))}{(1/\epsilon)^{d/2(\rho+d-1)}} .
$$
\end{thm}

\pf Following the same steps of the proof of Proposition \ref{thmlowergeo}, from formula (\ref{det}) to formula (\ref{J_m}), for any $\epsilon>0$, we attain the following lower bound
\begin{equation}\label{lower1}
   J_m:=  \frac{1}{2}\sum_{k=0}^{m} \tau_k^d \ln \left(a_k/\tau_k^d  \right) - \dim\mathcal{V}_m \ln\epsilon\leq \ln(\mathcal{C}(\epsilon, I_K )), \quad m=1,2,\ldots.
\end{equation}
For $m=1,2,\ldots$, we write 
\begin{equation}
I_m:=\frac{1}{2}\left(\ln a_0 + \sum_{k=1}^{m}\tau_k^d\ln a_k +\sum_{k=1}^{m}\tau_k^d \ln\left(1/\tau_k^d \right)\right)= \frac{1}{2}\sum_{k=0}^{m} \tau_k^d \ln \left(a_k/\tau_k^d \right),
\end{equation}
and if $a_0\geq 1$ we observe that, for $m=1,2,\ldots$, 
\begin{eqnarray*}
 I_m &\geq & \frac{1}{2}\left(\ln a_0 + (\dim\mathcal{V}_m -1)\ln c_2 -\rho(\dim\mathcal{V}_m -1)\ln m+\sum_{k=1}^{m}\tau_k^d \ln(1/\tau_k^d )\right).
\end{eqnarray*}
Since $\tau_k^d \leq \tau_m^d $ for $k=1,2\ldots,m$ we have 
$
\ln(1/\tau_k^d) \geq \ln(1/\tau_m^d ),
$
and, consequently 
\begin{eqnarray*}
 I_m &\geq & \frac{1}{2}\left[\ln a_0 + (\dim \mathcal{V}_m -1)\ln c_2 -\rho(\dim \mathcal{V}_m-1)\ln m+ (\dim \mathcal{V}_m-1)\ln(1/\tau_m^d )\right]
 \\ &=& \ln \sqrt{a_0}+ \frac{(\dim \mathcal{V}_m-1)}{2}\left[\ln c_2 -\rho\ln m+\ln(1/\tau_m^d )\right].
\end{eqnarray*}
From inequality (\ref{lower1}), we have 
\begin{equation}\label{Jm}
J_m=I_m- \dim\mathcal{V}_m \ln(\epsilon)\leq\ln(\mathcal{C}(\epsilon, I_K )), \quad m=1,2,\ldots.
\end{equation}
An application of Proposition \ref{dimVm} implies 
$$
\frac{\dim \mathcal{V}_m-1}{2}<\frac{\Gamma (\beta+1) m^d}{\Gamma(\alpha+2) \Gamma(\alpha+\beta +2)}, \quad \quad -\dim \mathcal{V}_m <-\frac{\Gamma(\beta+1) m^d}{2\,\Gamma(\alpha+2) \Gamma(\alpha+\beta+2)}, 
$$
and 
$$
\ln\left(\frac{1}{\tau_m^d}\right) \geq \ln \left(\frac{\Gamma(\alpha+1) \Gamma(\alpha+\beta+2)}{3\Gamma(\beta+1) m^{d-1}}\right),
$$
for $m=m_0, m_0+1, \ldots.$
Therefore,
\begin{equation}\label{base}
\ln \sqrt{a_0}+ \frac{\Gamma (\beta+1) m^d}{\Gamma(\alpha+2) \Gamma(\alpha+\beta +2) }\left[ \ln \left(\frac{c_2 \Gamma(\alpha+1) \Gamma(\alpha+\beta+2)}{3\Gamma(\beta+1)\, \epsilon^{1/2}}\right) -(\rho +d-1) \ln{m}\right]
\end{equation}
is less than $J_m$ for $m=m_0, m_0+1, \ldots$.

In order to proceed analogously to the final steps of the proof of the Theorem \ref{thmlowergeo}, we consider $\epsilon>0$ and $m=m(\epsilon)$ such that $1/(m+1)\leq \epsilon <1/m$, and
$$ 
\varphi_{\epsilon} (m):=  m^d\left[ \ln \left(\frac{c_2 \Gamma(\alpha+1) \Gamma(\alpha+\beta+2)}{3\Gamma(\beta+1)\, \epsilon^{1/2}}\right) -(\rho +d-1)\ln m\right].
$$
Due inequality (\ref{Jm}), it holds 
\begin{equation}\label{lower3}
\ln{\sqrt{a_0}}+\frac{\Gamma (\beta+1)}{\Gamma(\alpha+2) \Gamma(\alpha+\beta +2)} \,\varphi_\epsilon (m)< \ln(\mathcal{C}(\epsilon, I_K )), \quad m=m_0, m_0+1, \ldots,
\end{equation}
We observe that from inequality (\ref{lower3}), we will have 
\begin{equation}\label{lower5}
\ln{\sqrt{a_0}}+\frac{\Gamma (\beta+1)}{\Gamma(\alpha+2) \Gamma(\alpha+\beta +2)} \,\varphi_\epsilon (c)< \ln(\mathcal{C}(\epsilon, I_K )),
\end{equation}
for the critical point $c$ of $\varphi_{\epsilon}$.\ In this case the critical point is 
$$
c= \left(\frac{c_2 \Gamma(\alpha+1) \Gamma(\alpha+\beta+2)}{3 \Gamma(\beta+1)\, \epsilon^{1/2}}\right)^{1/(\rho+d-1)} e^{-1/d}, 
$$ and we have 
\begin{eqnarray*}
\varphi_\epsilon (c) = \frac{\rho+d-1}{d\,e}\left( \frac{c_2 \Gamma(\alpha+1) \Gamma(\alpha+\beta+2)}{3\Gamma(\beta+1)} \right)^{d/(\rho+d-1)} \left(\frac{1}{\epsilon}\right)^{d/2(\rho+d-1)},
\end{eqnarray*}
and the estimate follows. \eop

We close the paper by presenting an example of kernel fulfilling the assumptions of the theorems of this section. In what follows, for $\alpha, \beta >-1$ will write 
$$
c(\alpha,\beta):=\frac{\Gamma(\beta+1) }{\Gamma(\alpha+1)\Gamma(\alpha+\beta+2)}.
$$
We observe that $P_0^{\alpha,\beta}(s)=P_0^{\alpha,\beta}(1)=1$ implies that $J_0^{\alpha,\beta}(s)=1$. Then, a kernel $K$ that is an isotropic positive definite kernel given by the series expansion (\ref{WKlegendre}) with $a_0=1$, is written as
$$
K(x, y)=1+\sum_{k=1}^{\infty} a_k J_k^{(\alpha, \beta)}(\cos (d(x, y))), \quad x,y \in \mathbb{M}^ d.    
$$

\begin{ex}\label{exharmonicasymp} Let $\gamma>0$ and $K$ be the isotropic positive definite kernel given by the series expansion
\begin{equation}
K(x, y)=1+\sum_{k=1}^{\infty} k^{-d-\gamma} J_k^{(\alpha, \beta)}(\cos (d(x, y))), \quad x,y \in \mathbb{M}^ d,    
\end{equation} Then,
\begin{eqnarray*}
    \limsup_{\epsilon\to 0^+}\frac{\ln (\mathcal{C}(\epsilon, I_K ))}{(1/\epsilon)^{2d/(d+\gamma-1)}\ln(1/\epsilon)} \leq \frac{c(\alpha,\beta)}{(\alpha+1)}\frac{4^{d/(\gamma+d-1)}}{(\gamma+d-1)^{d/(\gamma+d-1)}}.
\end{eqnarray*}
and
\begin{eqnarray*}
    \frac{(\gamma+2d-1)}{e d [3c(\alpha,\beta)]^{d/(\gamma+2d-1)}}\frac{c(\alpha,\beta)}{(\alpha+1)} \leq \liminf_{\epsilon\to 0}\frac{\ln(\mathcal{C}(\epsilon, I_K ))}{(1/\epsilon)^{d/2(2d+\gamma-1)}},
\end{eqnarray*}
In particular, if $\mathbb{M}^d=\mathbb{S}^d$, then
\begin{eqnarray*}
    \limsup_{\epsilon\to 0^+}\frac{d! \ln (\mathcal{C} (\epsilon, I_K))}{2(1/\epsilon)^{2d/(d+\gamma-1)} \ln (1/\epsilon)} \leq \left(\frac{4}{\gamma+d-1}\right)^{d/(\gamma+d-1)},
\end{eqnarray*}
and
\begin{eqnarray*}
    \frac{(\gamma+2d-1)}{e d }\left(\frac{(d-1)!}{3}\right)^{d/(2d+\gamma-1)} & \leq & \liminf_{\epsilon\to 0}\frac{d! \ln(\mathcal{C}(\epsilon, I_K ))}{2 (1/\epsilon)^{d/2(2d+\gamma-1)}}.
\end{eqnarray*}

In fact, the first inequality is a an direct application of Theorem \ref{thmupperhar} and the second stated inequality follows from the direct application of Theorem \ref{thmlowerhar}, for $\rho = \gamma +d$. For the particular case $\mathbb{M}^d=\mathbb{S}^d$, we have
$$
c(\alpha,\beta) = c((d-2)/2,(d-2)/2) =\frac{1}{(d-1)!} \quad \mbox{and}\quad \frac{c(\alpha,\beta)}{(\alpha+1)} =\frac{2}{d!},  
$$
and the estimate follows.
\end{ex}

An interesting theoretical question related to bound covering numbers is to decide if the asymptotic limit $\lim_{\epsilon\to 0^+}\ln(\mathcal{C}(\epsilon, I_K ))/\varphi(\epsilon)$, exists for some function $\varphi: (0,\infty)\rightarrow\mathbb{R}$. In the affirmative case, the asymptotic limit deliveries us the best bounds for the numbers, and reciprocally precise bounds for the covering numbers could deliver the existence of the asymptotic limit. The estimates we present in Theorem \ref{thmlowergeo} and in Theorem \ref{thmuppergeo} are accurate but we are not able to ensure that the asymptotic limit exists from the bounds even on the spherical setting (Subsection \ref{subsec-spherical-gaussian}). The theme is explored in \cite{kuhn22-MR4474650} for approximation numbers of infinitely many times differentiable functions on $\mathbb{R}^d$ satisfying a size condition on the order of the derivatives on the compacts. 

\vspace{0.25cm}

\textbf{Acknowledgements.} This research is supported by the São Paulo Research Foundation (FAPESP) grant no. 2022/11032-0 and was financed in part by the Higher Education and Science Committee of the Republic of Armenia grant no. 23RL-1A027.


\bibliographystyle{plain}
\bibliography{references}

\bigskip

\noindent {\textsc{Karina Gonzalez}}\\
Faculty of Mathematics and Mechanics\\
Yerevan State University\\
0025 Yerevan, Republic of Armenia\\
\textsc{Email}: ngonzalezkarina@gmail.com

\medskip

\noindent \textsc{Thaís Jord\~{a}o}\\
Departament of Mathematics\\
Universidade de S\~{a}o Paulo\\
13566-590 S\~{a}o Paulo, Brazil\\
\textsc{Email}: tjordao@icmc.usp.br

\bigskip

\noindent{\textsc{Keywords}:} covering numbers; positive definite kernel, Reproducing Hilbert Space; isotropic kernel 

\medskip

\noindent{\textsc{MSC}:} 47B06; 42C10; 46E22; 41A60

\end{document}